\documentclass{amsart}

\usepackage{amsmath,amsfonts,latexsym,amssymb}

\usepackage{ifthen}

\ifthenelse{\pdfoutput > 0}%
{\usepackage[pdftex]{graphicx,color}}%
{\usepackage[dvips]{graphicx,color}}

\newcommand{\ra}{{\longrightarrow}}
\newcommand{\pg}{{\mathbb{P}(G)}}
\newcommand{\qz}{{\mathbb{Q}/\mathbb{Z}}}

\DeclareMathOperator{\Hom}{Hom}
\DeclareMathOperator{\hmg}{Hmg}
\DeclareMathOperator{\coc}{Coc}

\newtheorem{theorem}{Theorem}[section]
\newtheorem{definition}[theorem]{Definition}
\newtheorem{proposition}[theorem]{Proposition}
\newtheorem{corollary}[theorem]{Corollary}
\newtheorem{lemma}[theorem]{Lemma}

\date{June 2023}

\title{Homogeneous Functions and Algebraic $K$--theory}

\author{R. Keith Dennis and Reinhard C. Laubenbacher}

\address{R. Keith Dennis, Department of Mathematics, Cornell University, Ithaca, NY 14853}

\email{dennis@math.cornell.edu}

\address{Reinhard C. Laubenbacher, Department of Medicine, University of Florida, Gainesville, FL
32611}

\email{Reinhard.Laubenbacher@medicine.ufl.edu}



\begin{document}

\maketitle

%

\begin{abstract}
In this paper we develop the theory of homogeneous functions between finite abelian groups.  Here, a
function $f:G\longrightarrow H$ between finite abelian groups is homogeneous of degree $d$ if
$f(nx)=n^df(x)$ for all $x\in G$ and all $n$ which are relatively prime to the order of $x$.  We
show that the group of homogeneous functions of degree one from a group $G$ of odd order to $\qz$
maps onto $SK_1(\mathbb{Z}[G])$, generalizing a result of R. Oliver for $p$-groups.
\end{abstract}

\section{Introduction}

\bigskip

[Remark by first author:\newline

This paper was written December 12, 1994.  It was not published due to the procrastination of the
first author alone.  As the ideas have had applications \cite{KL} and should have additional ones, it is being
made available now.]




\bigskip

The concept of homogeneous function occured to the first author around 1975 as being the correct
framework for studying $SK_1$ of group rings.  It did not seem to inspire much interest from his
coauthors of \cite{ADS} or \cite{ADOS}, and the idea has now languished for almost twenty years with
nothing more than some scraps of paper outlining a few theorems.  Most recently the work of Yap
\cite{Ya} has rekindled interest in this project, leading to the present paper.  That homogenenous
functions should play a role here is (or at least should have been) evident in the proofs given in
\cite{De} for the classical situation of a vector space.  What was missing was the correct
definition, why it should be correct, and the courage to investigate the consequences.  It was not
until the theorem of Oliver \cite{ADOS} appeared that it was clear that this was really the correct
approach, although there is almost no hint of these ideas in that paper.  The correct description of
$SK_1$ of the integral group ring of a finite abelian group of odd prime power order first appears
in the thesis of Yap \cite{Ya}, along with the idea of the generalized transfer in this context.
The main ingredient that has been lacking until this point is the development of these ideas in a
categorical context to exploit the abundance of homogeneous functions for computations.  This paper
is meant to serve that purpose.  It seems reasonable to expect that homogenous functions of higher
degree should play a role in the higher $K$-theory of integral group rings of finite abelian groups.

The contents of the paper are as follows.  In Section 2 we develop the theory of homogeneous
functions, including explicit calculations; in Section 3 we define the transfer homomorphism.
Section 4 contains the relationship between homogeneous functions of degree one and $SK_1$ of group
rings of finite abelian groups of odd order.  The last section consists of a list of questions and
observations.

\section{Homogeneous Functions}

In this section we will work in the category of finite abelian groups, unless noted otherwise.

\begin{definition}
Let $G,\ H$ be finite abelian groups, and let $d$ be an integer.  If $d \geq 0$, then a function
$f:G \ra H$ is {\it homogeneous of degree} $d$ if $f(nx)=n^df(x)$ for all integers $n$ such that
$(n,o(x))=1$.  If $d < 0$, then $f$ is homogeneous of degree $d$ if $n^{-d}f(nx)=f(x)$ for all $n$
such that $(n,o(x))=1$.

Let $\hmg ^d(G,H)$ denote the set of homogeneous functions from $G$ to $H$ of degree $d$.
\end{definition}

Note that the set $\hmg ^d(G,H)$ is a finite abelian group under pointwise addition.  Viewed as a
functor from the category of finite abelian groups to itself $\hmg ^d(-,-)$ is covariant in the
second variable and contravariant in the first.

\begin{definition}
Define $\hmg ^d(G)=\hmg ^d(G,\qz)$ if $d \neq 0$, and $\hmg ^0(G)=\hmg ^0(G,\mathbb{Z})$.
\end{definition}

Note that for each positive integer $m$, the group $\qz $ contains exactly one cyclic group of order
$m$.  (One can think of $\qz $ as the direct limit of all finite cyclic groups with respect to
inclusion maps.)  Since any homogeneous function on $G$ takes values inside a finite cyclic
subgroup, the group $\hmg ^d(G), d \neq 0$ is isomorphic to $\hmg ^d(G,\mathbb{Z}/m)$ for all $m$ which
are divisible by the exponent of $G$.

We now give a different description of homogeneous functions from a more homological point of view,
which is also more amenable to calculation.

\begin{definition}
Let $G$ be a finite abelian group (written additively), and let $d$ be an integer. Let $[G]$ denote
the free abelian group on $G$ with generators $[x]$ for each $x \in G$.  Let $D$ be the subgroup of
$[G]$ generated by all elements of the form
\begin{equation*}
[nx]-n^d[x] \text{  for } x \in G,\ n \in \mathbb{Z},\ (n,o(x))=1,
\end{equation*}
if $d \geq 0$, and
\begin{equation*}
n^{-d}[nx]-[x],
\end{equation*}
if $d < 0$.  Now define
\begin{equation*}
G[d]:= [G]/D.
\end{equation*}
\end{definition}

By abuse of notation we will frequently denote elements of $G[d]$ by $[x]$ as well.

\begin{proposition}
For all $d\geq 0$, there is a canonical isomorphism
\begin{equation*}
\hmg ^d(G,H) \cong \Hom (G[d],H).
\end{equation*}
\end{proposition}

\begin{proof}{Proof}
Let $f\in \hmg ^d(G,H)$.  Define a homomorphism $\bar f:[G]\longrightarrow H$ by $\bar f([x])=f(x)$.
First let $d\geq 0$.  If $[nx]-n^d[x]$ is a generator of the subgroup $D$, then
\begin{equation*}
\bar f([nx]-n^d[x])=f(nx)-n^df(x)=0
\end{equation*}
since $f$ is homogeneous of degree $d$.  Thus, $f$ induces a homomorphism $\tilde f:G[d]
\longrightarrow H$.  It is easy to check that the map $f \mapsto \tilde f$ is a group homomorphism.

Now let $h:G[d]\longrightarrow H$ be a homomorphism.  Define $\tilde h \in \hmg ^d(G,H)$ by $\tilde
h(x)=h([x])$.  If $(n,o(x))=1$, then
\begin{equation*}
\tilde h(nx)=h([nx])=h(n^d[x])=n^dh([x])=n^d\tilde h(x).
\end{equation*}
Thus, $\tilde h \in \hmg ^d(G,H)$.  It is also straightforward to check that $h\mapsto \tilde h$ is
a homomorphism, and the two mappings are inverse to each other.

A similar proof works for $d<0$.  Finally, it is routine to verify that this isomorphism is
canonical in the sense that it induces an equivalence of functors $\hmg ^d(-,-)\cong \Hom (-[d],-)$.
This completes the proof of the proposition.
\end{proof}

\begin{corollary}
The group $\hmg (G)$ satisfies the following universal property.  Every $f\in \hmg ^d(G,H)$ can be
factored uniquely as the composition of the homogeneous function $G\ra G[d]$ given by $g \mapsto
[g]$ and a homomorphism $G[d] \ra H$.
\end{corollary}

Note that, if $d\neq 0$, then the function $G\ra G[d]$ is one-to-one, which immediately implies
surjectivity of the induced homomorphism $\Hom (G[d], H) \ra \hmg ^d (G,H)$ in Proposition 2.4.

\begin{corollary}
There are isomorphisms
\begin{itemize}
\item
$\hmg ^d(G,H_1\oplus H_2) \cong \hmg ^d(G,H_1)\oplus \hmg ^d(G,H_2)$;
\item
$\hmg ^d(G,H) \cong \hmg ^d(G)\otimes H$.
\end{itemize}
The first isomorphism is canonical.
\end{corollary}

\begin{proof}{Proof}
The first isomorphism is an immediate consequence of Proposition 2.4.  It reduces the computation of
$\hmg ^d(G,H)\cong \Hom (G[d],H)$ to the case where $H$ is cyclic.  Choose an embedding $i$ of $H$
into $\qz $.  Now define a homomorphism
\begin{equation*}
\Hom (G[d],\qz )\otimes H \ra \Hom (G[d], H)
\end{equation*}
by $g\otimes h \mapsto f$, where $f(x)=i^{-1}(g(x)\otimes i(h))$.  It is straightforward to see that
this map is an isomorphism, but it is non-canonical, since it depends on the choice of $i$.
\end{proof}

\begin{corollary}
There are canonical isomorphisms
\begin{itemize}
\item
$\hmg ^0(G)=\hmg ^0(G,\mathbb{Z}) \cong \Hom (G[0], \mathbb{Z})$;
\item
$\hmg ^d(G)=\hmg ^d(G, \qz)\cong \Hom (G[d],\qz)$
for $d\neq 0$.
\end{itemize}
\end{corollary}

\begin{proof}{Proof}
To see the first assertion, observe that a homogeneous function of degree zero on $G$ is a function
which is constant on the set of generators for each cyclic subgroup of $G$.  Thus, such a function
is given by arbitrarily assigning values to the elements of $G$, subject to this restriction.  But
it is straightforward to see that $G[0]$ is isomorphic to a free abelian group of rank equal to the
number of cyclic subgroups of $G$.  The assertion now follows.

The second assertion follows from the remark after Definition 2.2.
\end{proof}

Hence, for all $d\neq 0$ we have a (non-canonical) isomorphism $\hmg ^d(G) \cong G[d]$, since
$\widehat {G[d]}=\Hom(G[d],\qz) \cong G[d]$, non-canonically.

\smallskip
We now compute the structure of $G[d]$.

\begin{definition}
Let $d$ be an integer, and let $k$ be a positive integer.  If $d>0$, then define
\begin{equation*}
o_d(k)=gcd\{u^d-1\mid u\in \mathbb{Z}, u\equiv 1 \pmod k\}.
\end{equation*}
If $d<0$, then define $o_d(k)=gcd\{u^{-d}-1\}$, and if $d=0$, then define $u_0(k)=0$.

Let $G$ be a finite abelian group.  For an element $x\in G$ we define $o_d(x) = o_d(o(x))$.  One
easily sees that $o_1(x) = o(x)$.  In the future we will refer to the numbers $o_d(x)$ as the {\it
higher orders} of $x$.
\end{definition}

\begin{proposition}
Let $R(G)$ be a set of representatives of generators for all cyclic subgroups of $G$.  For all
integers $d$, there is a (non-canonical) isomorphism
\begin{equation*}
G[d] \cong \bigoplus _{x \in R(G)}\mathbb{Z} /\langle o_d(x)\rangle,
\end{equation*}
where $o_d(x)=o_d(o(x))$.
\end{proposition}

\begin{proof}{Proof}
First assume that $d>0$.  Define a mapping
\begin{equation*}
\varphi :\bigoplus _{x\in R(G)}\mathbb{Z}/\langle o_d(x)\rangle
          \longrightarrow G[d]
\end{equation*}
by $(1+\langle o_d(x)\rangle ) \mapsto [x]$.  To see that this is well defined observe that, if
$u\equiv 1 \pmod{o(x)}$, then
\begin{equation*}
(u^d-1)[x]=u^d[x]-[x]=[ux]-[x]=[x]-[x]=0
\end{equation*}
in $G[d]$.  Therefore $o_d(x)=o([x])$
in $G[d]$ divides $u^d-1$ for all such $u$,
hence it divides their greatest common divisor $o_d(x)$.
Thus $o_d(x)[x]=0$.

To define the inverse mapping, let
\begin{equation*}
\psi :[G] \longrightarrow \bigoplus _{x\in R(G)}\mathbb{Z}/\langle o_d(x)
 \rangle
\end{equation*}
be defined as follows.  Let $g\in G$ and let $x\in R(G)$ be the unique element such that $\langle
g\rangle =\langle x \rangle$, that is, there exists an integer $n$, relatively prime to $o(x)$ such
that $g=nx$.  Define $\psi ([g])=n^d+\langle o_d(x)\rangle$.  Now, let $[mg]-m^d[g]$ be a generator
of the subgroup $D$ of $[G]$.  Then $\langle mg \rangle =\langle g \rangle =\langle x \rangle$ for
some $x\in R(G)$.  If $g=nx$, then
\begin{equation*}
\psi ([mg]-m^d[g])=(nm)^d+\langle o_d(x)\rangle - m^dn^d+
\langle o_d(x) \rangle = 0.
\end{equation*}
Therefore, $\psi$ factors through $G[d]$, and it is straightforward to see that $\psi$ is inverse to
$\varphi$.  This completes the proof of the proposition for positive $d$.

A similar proof works for negative $d$.  If $d=0$, then the relations defining $G[0]$ are $[nx]-[x]$
for all $x\in G$ and $(n,o(x))=1$.  Thus, in $[G]$ we identify all cyclic summands coming from
generators of the same cyclic subgroup of $G$.  Hence, in $G[0]$, we obtain one copy of $\mathbb{Z}$ for
each cyclic subgroup of $G$.  This proves the proposition for $d=0$.
\end{proof}

\begin{corollary}
Let $d\neq 0$.  With notation as above, $o_d(x)$ is the order of $[x]\in G[d]$.
\end{corollary}

It remains to compute the higher orders of $x$ for $x\in R(G)$.

\begin{lemma}
Let $d, k$ be positive integers, and let $p$ be a prime.  Then
\begin{enumerate}
\item If $k|l$, then $o_d(k)|o_d(l)$; if $b|d$, then $o_b(k)|o_d(k)$;
\item
if $p|o_d(k)$, then $p|k$;
\item
we have $k|o_d(k)$;
\item if $p|\frac{o_d(k)}{k}$, then $p|d$;
\item if $p|k$, then $v_p(o_d(k)) \geq v_p(k)+v_p(d)$, and equality holds if $p$ is odd.  Here $v_p$
of an integer denotes the highest power of $p$ dividing it;
\item if $2|k$, then $v_2(o_d(k)) \leq v_2(k)+v_2(d)+1$, and equality holds if $k=2m$ with $m$ odd;
\item if $(k,k')=1$, then $o_d(kk')\leq o_d(k)o_d(k')$;
\item if $(d,d')=1$, then
\begin{equation*}
\frac{o_{dd'}(k)}{k}=\frac{o_d(k)}k\frac{o_{d'}(k)}k;
\end{equation*}
\item for $s\geq 1$ we have
\begin{equation*}
  o_{p^s}(k)=
  \begin{cases}
    k(p,k)^s & \text {if $p$ is odd,} \\
    k(2,k)^{s-1}(4,2+k)& \text{if $p=2$};
  \end{cases}
\end{equation*}
\item we have $lcm(o_d(k),o_{d'}(k))=o_{lcm(d,d')}(k)$ and $(o_d(k),o_{d'}(k))=o_{(d,d')}(k)$;
\item $o_d(k)=o_{-d}(k)$.
\end{enumerate}
\end{lemma}

\begin{proof}
1. If $k|l$ and $u\equiv 1 \pmod l$, then $u\equiv 1 \pmod k$.
The first assertion is now clear.

If $d=nb$ and $u\equiv 1 \pmod k$, then $u^n\equiv 1 \pmod k$.
Furthermore, $u^d-1=(u^n)^b-1$.  The second statement is now clear.

2. Suppose that $p$ does not divide $k$.  Then there is some $s$
such that $ps\equiv 1 \pmod k$.  But $o_d(k)|((ps)^d-1)$, therefore
$p$ cannot divide $o_d(k)$.

3. This assertion is clear since $k$ divides any number of the form
$(ak+1)^d-1$.

4. Write $(ak+1)^d-1=\sum _{i=0}^d \binom di (ak)^i -1 =akE$,
where
\begin{equation*}
E=d+\binom d2 ak+\binom d3 a^2k^2 + \cdots + \binom d{d-1}a^{d-2}k^{d-2}
  +a^{d-1}k^{d-1}.
\end{equation*}
If $p|\frac{o_d(k)}k$, then $p|(k+1)^d-1$, so that $p|E$.  Since
$p|k$ it follows that $p|d$.

5. Let $d=p^tb$ with $(b,p)=1$.  Then
\begin{equation*}
v_p(\binom {p^tb}ia^{i-1}k^{i-1})\geq t+\epsilon +(i-1)-\sum _i[\frac i{p^i}].
\end{equation*}
To see this, first observe that,
since $p|k$, we have that $p^{i-1}|k^{i-1}$.  This accounts for the
summand $i-1$ on the right-hand side of the inequality.

Furthermore, we have
\begin{equation*}
\frac{(p^tb)!}{(p^tb-i)!}=(p^tb-i+1)\cdots (p^tb).
\end{equation*}
Hence the right-hand side is divisible at least by $p^{t+\epsilon}$,
where $\epsilon =1$ if $i>p$ and $\epsilon =0$ otherwise.  Since
$v_p(i!)=\sum_{j\geq 1}[\frac i{p^j}]$ by the lemma below, we obtain the above
inequality.

Also, by the lemma below, if $i\geq 2$, then $\sum _j[\frac i{p^j}] \leq i-1$.  Thus, we have shown
that every term of $E$ in (4), except for the first, possibly the second (if $p=2=i$), and possibly
the last, is divisible by $p^{t+1}$.  For the last term, we see that $v_p(a^{d-1}k^{d-1} \geq d-1$,
since $p|k$.  And $d-1 =p^tb-1 \geq t+1$ unless $d=1$ or $p=d=2$.  Thus, in all cases, all terms are
divisible by $p^t$, and all except the first and possibly the second are divisible by $p^{t+1}$.

For $p$ odd, all terms except for the first are divisible by $p^{t+1}$, and the first term is
divisible by $p^t$ only.  Hence, the largest power of $p$ that divides $E$ is $p^t$, when $a=1$.
This proves (5).

6. Let $a=2=p$ in the previous discussion.  Then it follows that the exact power of $2$ dividing
$akE$ is $1+v_2(k)+v_2(d)$, since each term of $E$ except for the first one is divisible by
$2^{t+1}$.  This yields the inequality.

As before, to determine whether the exact power of $2$ dividing $akE$ is $2^t$ or $2^{t+1}$ we only
need to consider the first two terms of $E$.  If $k$ is divisible by 4, then the second term is
divisible by $2^{t+1}$.  Since the first term is only divisible by $2^t$, it follows that $E$ is
divisible exactly by $2^t$.  If, however, $k=2m$ with $m$ odd, then we have
\begin{equation*}
d + \binom d2 a\cdot 2m = d(1+(d-1)am).
\end{equation*}
If $t=0$, then $E=1$, hence $akE$ is even, since $k$ is even.  Now suppose that $t>0$.  If $a$ is
even, then $1+(d-1)am$ is odd, hence $E$ is divisible by $2^t$, and $akE$ is divisible by $2^{t+1}$.
If $a$ is odd, then $1+(d-1)am$ is even, so that $E$ is divisible by $2^{t+1}$.  Thus we obtain
equality in this case.

7.
Let $p$ be a prime divisor of $o_d(kk')$.  If $p$ is odd, then
\begin{equation*}
v_p(o_d(kk'))=v_p(kk')+v_p(d)=v_p(k)+v_p(k')+v_p(d)
\leq v_p(o_d(k))+v_p(o_d(k')).
\end{equation*}
If $p=2$, then
\begin{equation*}
v_2(o_d(kk'))\leq v_2(kk')+v_2(d)+1=v_2(k)+v_2(k')+v_2(d)+1,
\end{equation*}
and equality holds if $kk'=2m$ with $m$ odd.  The inequality now
follows since \break $v_2(o_d(kk')\geq v_2(kk')+v_2(d)$ and $k$ and $k'$
are relatively prime.

8.
If $p$ is a prime divisor of $o_d(k)$, then $p|k$.  If $p$ is odd,
then it follows from (5) that $v_p(o_d(k))=v_p(k)+v_p(d)$.  Thus,
$v_p(o_{dd'}(k))=kv_p(o_d(k))v_p(o_{d'}(k))$.  Similarly, if
$p=2$, then the same argument works, with an extra factor of 2,
if $k=2m$, $m$ odd.

9. This formula follows immediately from (5) and (6), upon
noting that $(4,2+k)$ is equal to 1 when $k$ is odd, equal to 2 if
$k$ is divisible by 4, and equal to 4 if $k$ is divisible by 2 but
not 4.

10. The formulas for least common multiple and greatest common divisor
follow immediately from (7) and (8).

11. Obvious.
\end{proof}

\begin{lemma}
Let $n$ be a positive integer and $p$ a prime.  Then
 $v_p(n!)=\sum _i [\frac i{p^i}]$.  Furthermore,  $\sum _j[\frac i{p^j}] \leq i-1$
and equality holds precisely when $p=2$ and $n=2^r$ for some non-negative
integer $r$.
\end{lemma}

\begin{proof}{Proof}
 The first assertion is straightforward to show.

To verify the second assertion,
let $p^n$ be the highest power of $p$ less than or equal
to $i$.  Then
\begin{equation*}
\begin{split}
[\frac ip]+[\frac i{p^2}] + \cdots +[\frac i{p^n}]& \leq \frac ip +
\cdots +\frac i{p^n} \\
&=\frac{p^{n-1}+p^{n-2}+\cdots +p+1}{p^n}\cdot i\\
&=(\frac 1{p-1}-\frac 1{(p-1)p^n})\cdot i \\
&< i.
\end{split}
\end{equation*}
This shows the inequality.  Now suppose that $\sum _j[\frac i{p^j}]=i-1$.
Then
\begin{equation*}
i-1 \leq (\frac 1{p-1} - \frac 1{(p-1)p^n})\cdot i,
\end{equation*}
which implies that $(p-1)(1-\frac 1i) \leq 1-\frac 1i$, so
that $p-1 \leq 1$.  Hence $p=2$.  Thus
\begin{equation*}
i-1 \leq \frac {2^n-1}{2^n} \cdot i = i-\frac 1{2^n}i,
\end{equation*}
which implies that $i\leq 2^n$.  Since $2^n$ is the largest power of
$2$ which is less than or equal to $i$, we have equality.  This completes
the proof of the lemma.
\end{proof}

\begin{theorem}
Let $G$ be a finite abelian group, and $d\neq 0$.
Then $\hmg ^d(G)$ has the following
canonical decomposition into Sylow subgroups:
\begin{equation*}
\hmg ^d(G) \cong \bigoplus _{p\mid |G|} \hmg ^d(G_p)^{q(G/G_p)},
\end{equation*}
where $G_p$ is the $p$-Sylow subgroup of $G$, and $q(G/G_p)$
is the number of cyclic subgroups of $G/G_p$.
\end{theorem}

\begin{proof}{Proof}
The proofs for positive and negative $d$ are similar.  Thus,
we assume that $d>0$.
In light of the canonical isomorphism
\begin{equation*}
\hmg ^d(G)\cong \Hom (G[d],\qz)
\end{equation*}
it is sufficient to show that we have a canonical decomposition
\begin{equation*}
G[d] \cong \bigoplus _{p\mid |G|}G_p[d]^{q(G/G_p)}.
\end{equation*}
For a cyclic group $C$ let $gen(C)$ be the set of generators of
$C$.  Now observe that
\begin{equation*}
[G]=\bigoplus_{ \substack{C \subset G \\ \text{ cyclic}}}
\left(\bigoplus _{x\in gen(C)}
\mathbb{Z}[x]\right)=\bigoplus _C[G]_C.
\end{equation*}
Furthermore, there is a similar decomposition $D=\bigoplus _CD_C$,
where $D_C$ is the subgroup of $D$ generated by the relations
$[nx]-n^d[x]$ with $x\in gen(C)$.  Then $D_C \subset [G]_C$
and there is a canonical
isomorphism
\begin{equation*}
G[d]\cong \bigoplus_{\substack{C\subset G \\\text { cyclic}}} [G]_C/D_C,
\end{equation*}
and $[G]_C/D_C$ is a finite cyclic group.  We get a similar decomposition
for
\begin{equation*}
G_p[d]\cong \bigoplus_{\substack{C\subset G_p \\\text{ cyclic}}} [G]_C/D_C.
\end{equation*}

Now consider $[G]_C/D_C$ for $C \nsubseteq G_p$.  Let $C\cong P\oplus Q$,
where $P$ is the $p$-Sylow subgroup of $C$.  We have that
\begin{equation*}
[G]_C = \bigoplus _{x\in gen(C)}\mathbb{Z}[x].
\end{equation*}
There is a one-to-one correspondence between generators of $C$
and pairs $(a,b)$ such that $a$ is a generator of $P$ and $b$ is
a generator of $Q$.
Now define a homomorphism
\begin{equation*}
\varphi :[G]_C/D_C \ra [G]_P/D_P \oplus [G]_Q/D_Q
\end{equation*}
by $[x] \mapsto ([a],[b])$ for every generator $x=a+b$ of $C$.  If
$n$ is relatively prime to $o(x)$, then it is also relatively prime
to $o(a)$ and $o(b)$, so that $\varphi $ is well-defined.  Since
every generator $([a],[b]$ of $[G]_P/D_P \oplus [G]_Q/D_Q$ is
in the image of $\varphi$, it is onto.  Now observe that
$o([x])=o([a])o([b])$ in $G[d]$.  But
\begin{equation*}
o([x])=o_d(o(x))=o_d(o(a)o(b)) \leq o_d(o(a))o_d(o(b))=o([a])o([b]),
\end{equation*}
by Proposition 2.11. (7)
This shows that the source and target of $\varphi $ have the same
order so that it is an isomorphism.
Thus the cyclic group
$[G]_C/D_C$ contains a subgroup canonically isomorphic to
$[G_P]/D_P$.  Fixing $Q$ and letting $P$ range over all cyclic
subgroups of $G_p$, we see that $G[d]$ contains a subgroup
canonically isomorphic
to $G_p[d]$ for every
cyclic subgroup $Q$ of $G/G_p$.  The theorem now follows.
\end{proof}

\begin{corollary} \cite{Ya}
Let $G$ be a finite abelian group.  Then there is a canonical
isomorphism of groups
\begin{equation*}
\hmg (G)=\hmg ^1(G)\cong \bigoplus_{\substack{C \subset G \\\text{ cyclic}}}
\widehat C,
\end{equation*}
where $\widehat C=\Hom (C,\qz)$ is the dual of $C$.
\end{corollary}

\begin{proof}{Proof}
We showed in the proof of the previous theorem that
\begin{equation*}
G[d]\cong \bigoplus_{\substack{C \subset G \\\text{ cyclic}}} [G]_C/D_C,
\end{equation*}
where
\begin{equation*}
[G]_C=\bigoplus _{x\in gen(C)}\mathbb{Z}[x].
\end{equation*}
If $d=1$, then it
is easy to see that there is a canonical isomorphism
$C\cong [G]_C/D_C$, given by $x \mapsto [x]$ for a generator $x$ of $C$.
The corollary now follows.
\end{proof}

We now study a subgroup of $\hmg (G)$ which will play an important
role in the relationship of $\hmg (G)$ to $K_1(\mathbb{Z}[G])$.

\begin{definition}
Let $K\subset G$ be a subgroup of a finite abelian group.  Then
$K$ is called {\it cocyclic} if $G/K$ is a cyclic group.

Let $\phi :K\longrightarrow \qz$ be a character of $K$.
Then $\phi$ induces a function $\phi _K \in \hmg ^1(G)=\hmg (G)$
which is defined by
\begin{equation*}
\phi _K(g)=\begin{cases} \phi (g)& \text{if $g\in K$}\\
                  0 &\text{if $g\notin K$}.
           \end{cases}
\end{equation*}
Call $\phi _K$ a {\it cocyclic function}.

Let $\coc (G)$ denote the subgroup of $\hmg (G)$ generated by all
cocyclic functions.
\end{definition}

\begin{theorem}
Let $G$ be a finite abelian group.  Then there is a canonical
isomorphism
\begin{equation*}
\hmg (G)/\coc (G) \cong
\bigoplus _{p\mid |G|}\left(\hmg (G_p)/\coc (G_p)\right)
^{q(G/G_p)},
\end{equation*}
where $q(G/G_p)$ denotes the number of cyclic subgroups of $G/G_p$.
\end{theorem}

\begin{proof}{Proof}
In light of Theorem 2.13 it is sufficient to show that
there is a canonical
decomposition into Sylow subgroups:
\begin{equation*}
\coc (G) \cong \bigoplus _{p\mid |G|} \coc (G_p)^{q(G/G_p)}.
\end{equation*}
Let $K$ be a cocyclic subgroup of $G$.
Under the canonical isomorphisms
\begin{equation*}
\hmg (G) \cong \Hom (G[1],\qz)
\text{  and  } G[1] \cong \bigoplus_{\substack{C\subset G \\ \text{ cyclic}}} C,
\end{equation*}
the function $\phi _K$ corresponds to the homomorphism
\begin{equation*}
\varphi _K: \bigoplus _C C \longrightarrow \qz
\end{equation*}
given by
\begin{equation*}
\varphi _K|_C(n[x]) = \begin{cases} n\phi (x)=\phi (nx)&
                             \text{if $C=\langle x \rangle \subset K$},\\
                             0&\text{otherwise}.
                      \end{cases}
\end{equation*}

Recall now the decomposition of $G[1]$ into Sylow subgroups:
\begin{equation*}
G[1]\cong \bigoplus _{p\mid |G|}G_p[1]^{q(G/G_p)}.
\end{equation*}
Let $p$ be a prime divisor of $|G|$ and let $Q$ be a cyclic
subgroup of $G/G_p$.  Let $P$ be a cyclic subgroup of $G_p$, so
that $P\oplus Q = C$ is a cyclic subgroup of $G$.  Conversely,
every cyclic subgroup of $G$ can be decomposed in this way.
If $Q$ is not contained in $K$, then $P\oplus Q$ is not a subgroup of $K$
either, so that $\varphi _K|_{P\oplus Q}=\varphi _K|_{P\oplus 0}=0$.
Thus, $\varphi _K$ is the zero function on the the copy of
$G_p[1]$ coming from $Q$ in the above decomposition of $G[1]$,
which is a cocyclic function on $G_p[1]$.

Now suppose that $Q \subset K$.  If $P$ is not contained in $K$,
then $\varphi _K|_{P\oplus Q}=0$ as before.  So assume that
$P\oplus Q \subset K$.  Let
$P=\langle x \rangle , Q=\langle y \rangle$,
with $o(x)=n, o(y)=m$.  Then $P\oplus Q=\langle x+y \rangle $, and
$\langle an(x+y)\rangle $ is the subgroup of $P\oplus Q$ canonically
isomorphic to $P$, where $a,b$ are integers such
that $an+bm=1$.  Then $an(x+y)=anx=(1-bm)x=x$, and
\begin{equation*}
\varphi _K(an[x+y])=an\phi _K (x+y)=\phi _K(an(x+y))=\phi _K(x).
\end{equation*}
Hence $\varphi _K|{\langle an[x+y]\rangle} =\phi _K|_P$.
In summary, when restricted to the summand $G_p[1]$ of $G[1]_p$
coming from $Q$,
the function $\varphi _K$ is zero on cyclic summands that are
not contained in $G_p\cap K$ and $\phi _K|_{G_p}=\phi _{K\cap G_p}$
otherwise.  Thus, $\varphi _K$ restricted to  any of the
summands $G_p[1]$ is a cocyclic function.  Therefore, every
cocyclic function $\varphi _K$ on $G[1]$ decomposes as
a sum of cocyclic functions
\begin{equation*}
\varphi _K = \sum _{p\mid |G|}\sum_{\substack{Q\subset G/G_p \\ \text{cyclic}}}
               \varphi _K|_{G_p[1]}.
\end{equation*}
Thus, we obtain a one-to-one homomorphism
\begin{equation*}
\Phi :\coc (G) \ra \bigoplus _{p\mid |G|} \coc (G_p)^{q(G/G_p)}.
\end{equation*}

To show that this homomorphism is onto
let $K_p$ be a cocyclic subgroup of $G_p$, and let
$\phi _{K_p}$ be one of the generators of the summand $\coc (G_p)$ belonging
to the cyclic subgroup $Q \subset G/G_p$.
That is, $\phi _{K_p}$ is induced from a character $\phi :K_p\ra \qz $.
Now let $K$ be the cocyclic subgroup $K_p\oplus G/G_p$ of $G$, and let
$\psi :K \ra \qz $ be the character defined by $\phi $ on the first
summand and by the character induced by $Q$ on the second summand.
Then it is straightforward to see that the cocyclic function $\psi _K$
maps to $\phi _{K_p}$ under $\Phi $.
This completes the proof of the theorem.
\end{proof}

\section{The Transfer}

In this section we define a transfer homomorphism
\begin{equation*}
T_t:\hmg^d(G,H) \longrightarrow \hmg^d(G',H)
\end{equation*}
induced by a homogeneous function $t:G \longrightarrow G'$ of degree $d$.  This generalizes the
transfer induced by a group homomorphism, defined in \cite{Ya}.  Let $G,G'$ be finite abelian groups
and assume that $G$ has odd order.  Let $t\in Hmg^d(G,G')$.  By Proposition 2.4, it is enough to
show that $t$ induces a homomorphism
\begin{equation*}
T_t:\Hom(G[d],H)\longrightarrow \Hom(G'[d],H).
\end{equation*}

Observe first, that $t$ induces a homomorphism $[G]\longrightarrow [G']$,
which we will again denote by $t$.  This homomorphism has the
property that the intersection of the image of $t$ with each cyclic
factor of $[G']$ is either equal to zero or is equal to the cyclic
factor.  Now decompose $G'[d]$ as in the proof of Theorem 2.13:
\begin{equation*}
G'[d] \cong \bigoplus_{\substack{C\subset G \\ C \text{ cyclic}}}
[G']_C/D'_C =\bigoplus _C A_C.
\end{equation*}
Let $C$ be a cyclic subgroup of $G'$.  If $x$ is a generator of $C$
which is in the image of $t$, then the induced homomorphism
\begin{equation*}
\tilde t:G[d] \longrightarrow G'[d]
\end{equation*}
maps onto $[G']_C/D'_C=A_C$, since it is generated by the residue class
of $[x]$.
If no generator of $C$ is in the image of $t$, on the other hand,
then the intersection of the image of $t$ with $[G']_C$ is zero,
hence $im(\tilde t)\cap A_C=0$.

Now we are ready to define the transfer $T_t$.

\begin{definition}
Let $f\in \Hom(G[d],H)$.
In order to define $T_t(f)\in \Hom(G'[d],H)$, it is sufficient to define
it on each cyclic factor $A_C$ of $G'[d]$.  Let $x\in A_C$.
Define
\begin{equation*}
T_t(f)(x)=\sum_{\substack{y\in G[d] \\ \tilde t(y)=x}} f(y).
\end{equation*}
Here the empty sum is to be interpreted as zero.
\end{definition}

We need to prove
that $T_t(f)$ is a homomorphism.

If $im(\tilde t)\cap A_C=0$, then $T_t(f)=0$ on $A_C$, since
the order of $G$, hence the order of $G[d]$, is odd, so
that $\sum {y\in G[d]} f(y) =0$.  Now suppose that $\tilde t$
maps onto $A_C$, and let $x,x'\in A_C$.
First observe that, if $y_0\in G[d]$ with $\tilde t(y_0)=x$, then
$\tilde t^{-1}(x)=\tilde t^{-1}(0)+y_0$; similarly for $x'$ and
$x+x'$.  Therefore,
\begin{equation*}
\sum {\tilde t(y)=x}f(y) =\sum _{z\in \tilde t^{-1}(0)}f(z+y_0)
  =af(y_0)+\sum _{z \in \tilde t^{-1}(0)}f(z) =af(y_0),
\end{equation*}
where $a=|\tilde t^{-1}(0)|$.  The last sum is zero, since it is
the sum of all elements in an abelian group of odd order.  Thus,
if $\tilde t(y'_0)=x'$ and $\tilde t(y_0+y'_0)=x+x'$,
then
\begin{equation*}
T_t(f)(x+x')=af(y_0+y'_0)=af(y_0)+af(y'_0)=T_t(f)(x)+T_t(f)(x').
\end{equation*}
Hence $T_t(f)\in \Hom(G'[d],H)$.

\begin{proposition}
Let $t:G\longrightarrow G'$ be homogeneous of degree $d$, and assume
that $|G|$ is odd.  Let
$\tilde t:\Hom(G[d],H)\longrightarrow \Hom(G'[d],H)$
and $t^*:\Hom^d(G',H)\longrightarrow \Hom^d(G,H)$ be induced by $t$
as before.
\begin{itemize}
\item
If $t$ is one-to-one, then $T_t$ is one-to-one;
\item
the composition $t^*\circ T_t$ is multiplication by $|ker(\tilde t)|$;
\item
the composition $T_t\circ t^*$ is given by
\begin{equation*}
(T_t\circ t^*)(g)(x)=T_t(g\circ \tilde t)(x)=
\begin{cases} 0&\text{if
                   $x\notin im(\tilde t)$} \\
|ker (\tilde t)|g(x)& \text{otherwise}.
\end{cases}
\end{equation*}
\end{itemize}
\end{proposition}

\begin{proof}{Proof}
All three assertions follow directly from the observation made
earlier, that \break $T_t(f)(x)=|ker(\tilde t)|f(y)$ for any $y\in G[d]$
such that $\tilde t(y)=x$, if $x\in im(\tilde t)$ and zero otherwise.
\end{proof}

\section{Homogeneous functions and $SK_1$ of group rings}

Following is the main result of this section, which generalizes
a result of R. Oliver (Theorem 4.2 below) for $p$-groups.

\begin{theorem}
Let $G$ be a finite abelian group of odd order.  Then there is
an isomorphism
\begin{equation*}
\hmg (G)/\coc (G) \cong \hmg (\widehat G)/\coc (\widehat G)\cong SK_1(\mathbb{Z}[G]).
\end{equation*}
\end{theorem}

In order to prove this theorem we recall some results from
\cite{ADS} and \cite{ADOS}.
The first result reduces the computation of $SK_1(\mathbb{Z}[G])$
to the case of a $p$-group.  It inspired the decomposition
in Theorem 2.13 above.

\begin{theorem} \cite[Theorem 3.11]{ADS}
Let $G_p$ be the $p$-Sylow subgroup of $G$.  Then
\begin{equation*}
SK_1(\mathbb{Z}[G]) \cong \prod _{p\mid |G|}
   \left(SK_1(\mathbb{Z}[G_p])\right)^{q(G/G_p)}.
\end{equation*}
\end{theorem}

Combined with Theorem 2.16, this result reduces the proof of Theorem 4.1 to the case of a $p$-group.
Let $\chi $ be an irreducible character of $G$.  Then $\chi$ corresponds to a simple factor of $\mathbb{Q}[G]$, which is isomorphic to a cyclotomic extension $\mathbb{Q}(\chi )$ of $\mathbb{Q}$, obtained by adding
a primitive $|im(\chi )|$-th root of unity to $\mathbb{Q}$.  Furthermore, $\chi$ induces a surjective
homomorphism from $\mathbb{Q}[G]$ to $\mathbb{Q}[\chi ]$.  Denote by $\mathbb{Z}[\chi ]$ the image of $\mathbb{Z}[G]$ in $\mathbb{Q}(\chi )$ under the composition of this homomorphism with the inclusion.  Let
$R=R(G)$ denote a set of representatives of all irreducible characters of $G$, and let $R_0$ denote
the subset of all nontrivial irreducible characters.

\begin{theorem} \cite[Theorem 2.10]{ADS}
Let $G$ be a finite abelian $p$-group of order $p^l$.  Then for
any $k\geq l$ there is an exact sequence
\begin{equation*}
K_2(\mathbb{Z}[G]/(p^k)) \ra \prod _{\chi \in R_0} K_2(\mathbb{Z}[\chi ]/(p^k))
  \ra SK_1(\mathbb{Z}[G]) \ra 1\, .
\end{equation*}
\end{theorem}

In \cite[Sect. 1]{ADOS} it was shown that the middle term in
this exact sequence is isomorphic to $\prod _{\chi \in R_0} im(\chi )$.
Define a homomorphism
\begin{equation*}
\varphi :\hmg (\widehat G) \ra \prod _{\chi \in R_0} im(\chi )
\end{equation*}
by $\varphi (f) _{\chi }= f(\chi )$.  It is clear that $\varphi $
is one-to-one.  Furthermore, it maps onto each factor of the product.
To see this, let $\chi $ be an index in the product and let $x$
be a generator of $im(\chi )$.  Then define
$f\in \hmg (\widehat G)$ by
\begin{equation*}
f(\eta )= \begin{cases} nx &\text{if $\eta =n\chi $,}\\
                 0  &\text{otherwise.}
          \end{cases}
\end{equation*}
Then $\varphi (f) = x$.  Therefore $\varphi $ is
an isomorphism.  This completes the proof of Theorem 4.1.

\bigskip
For the case of a $p$-group this result is used in \cite{Ya} to
obtain information about $SK_1$ of an elementary abelian $p$-group
and $\mathbb{Z}/p^e\oplus \mathbb{Z}/p^e$.

\section{Questions and Speculations}

We end this paper with a list of questions and observations which we plan to pursue in a subsequent
paper.  It is our hope that homogeneous functions of higher degree are related to the higher
$K$-theory of integral group rings of finite groups and might help answer some of the many questions
remaining even in dimension one.  \smallskip \noindent 1. Do there exist maps
\begin{equation*}
\hmg ^d(\widehat G) \ra K_d(\mathbb{Z}[G])
\end{equation*}
or
\begin{equation*}
\hmg ^d(\widehat G) \ra K_d(\mathbb{Q}[G])?
\end{equation*}
For $d=0$, there is of course an abstract isomorphism $\hmg ^0(G)\cong K_0(\mathbb{Q}[G])$.  The
situation in dimension one suggests that in higher dimensions one most likely needs to first develop
the appropriate number theoretic machinery.  To start with, one might try to make the isomorphism in
Theorem 4.1 explicit, that is, associate to a homogeneous function in $\hmg (\widehat G)$ an
explicit matrix in $SK_1(\mathbb{Z}[G])$.

\smallskip
\noindent
2. Do the groups $\hmg^d(G)$ assemble to a cohomology theory?  Do
they form a graded ring?

\smallskip
\noindent
3. Is there a theory of homogeneous functions for non-abelian groups?
A reasonable starting point would be class functions.  Can such a
theory be used to describe $SK_1(\mathbb{Z}[G])$ (or the class group
$Cl(\mathbb{Z}[G])$ of Oliver)?

\smallskip
\noindent
4. In what sense do the groups $\hmg ^d(G)$ reflect the rational
representation theory of the group $G$?  For instance, are they
related to the Burnside ring of $G$?

\smallskip
\noindent
5. Can the results in Section V be generalized to coefficient rings
other than $\mathbb{Z}$?  In \cite{ADS} and \cite{ADOS} it was shown
that in many cases (for instance for rings of algebraic integers
in totally real extensions of $\mathbb{Q}$ in which all prime divisors
of $|G|$ are unramified) the group $SK_1(\mathbb{Z}[G])$ does not
depend on the coefficient ring.

\smallskip
\noindent
6. Let $U$ denote the group of multiplicative units in the ring
$\widehat {\mathbb{Z}}$ defined to be the inverse limit of the rings
$\mathbb{Z}/m$ taken with respect to the directed set given by divisibility
of integers.  Let $G$ be a finite abelian group of order $m$.  Then
$(\mathbb{Z}/m)^*$ acts on $G$ in a natural way via multiplication,
making $G$ into a $(\mathbb{Z}/m)^*$-module.  Since there is a
unique surjective ring homomorphism $U \ra (\mathbb{Z}/m)^*$ for
every integer $m$, it follows that $G$ is a $U$-module (and
a homogeneous function on $G$ of degree one is nothing but
a $U$-module homomorphism).

Let $\pg$ denote the set of orbits of $G$ under the action of $U$.
One might think of $\pg$ as a {\it projective space on} $G$.
There is a one-to-one correspondence between the elements of
$\pg$ and the cyclic subgroups of $G$, so that the
number of ``points'' in $\pg$ is equal to the number of
irreducible rational representations of $G$.

Is there a ``projective geometry'' on $\pg$?  Can one use
it to develop a geometric/combinatorial method to
compute $SK_1(\mathbb{Z}[G])$ for finite abelian groups?  For instance,
one might study filtrations on $\mathbb{Z}[G]$ arising from cyclic
decompositions of $G$.

\smallskip
\noindent
7. Can one extend the results in this paper to groups of even order?
Oliver's result in this case gives a description of $SK_1(\mathbb{Z}[G])$
as a subquotient of $\hmg (\widehat G)$ (in our terminology).
For a small number of computations the description is isomorphic
to that given for odd primes.

\smallskip
\noindent
8. Can one use the theory developed in this paper to prove the
conjecture in \cite{ADOS} that for all odd primes $p$, and a
fixed list of invariants for $G$, the invariants for $SK_1(\mathbb{Z}[G])$
are given by rational polynomial functions in $p$ which are
{\it independent} of the prime $p$?

\end{document}